\documentclass[12pt,a4paper]{article}

\usepackage{amssymb}
\usepackage{amsmath}

\newtheorem{Proposition}{Proposition}[section]
\newtheorem{Theorem}[Proposition]{Theorem}
\newtheorem{Corollary}[Proposition]{Corollary}

\newcommand{\C}{{\mathbb C}}
\newcommand{\N}{\mathbb N}
\newcommand{\Z}{\mathbb Z}

\newcommand{\Qp}{\mathbb Q_p} 
\newcommand{\Ql}{\mathbb Q_\ell} 
\newcommand{\Zl}{\mathbb Z_\ell} 
\newcommand{\Cp}{\mathbb C_p} 
\newcommand{\Zln}{\Z/\ell^n\Z} 
\newcommand{\Zn}{\Z/n\Z} 

\newcommand{\Q}{\mathbb Q} 
\newcommand {\F}{\mbox{${\cal F}$}} 
\newcommand {\Fl}{\mbox{${\cal F} \!\!\!\!\!\!\!\; {\cal J}$}}

\newcommand{\proof}{\noindent {\bf Proof: }}
\newcommand {\qed}{\hfill $\square$}

\begin{document}

\title{The continuous cohomology of period domains over local fields}

\author{Sascha Orlik}

\date{}

\maketitle

\bigskip
\noindent {\bf Abstract.}
In this paper we consider period domains over local fields for quasi-split reductive groups. We compute the continuous $\ell$-adic cohomology with compact support of them in the case of a basic isocrystal. 
This paper is a continuation of \cite{O2} where we considered the \'etale cohomology with torsion coefficients of these spaces. 


\bigskip
\section{Introduction}

This paper deals with the computation of the continuous $\ell$-adic 
cohomology with compact support of period domains over local fields. We consider the case where the underlying reductive group is quasi-split  and the isocrystal is basic. In \cite{O2} we computed the \'etale cohomology with torsion coefficients of them.  
If we pass to the usual $\ell$-adic cohomology defined by the projective limit of the cohomology groups with
values in $\Zln ,$ we get continuous generalized Steinberg representations in Banach spaces. This was first pointed out by R. Huber in the Drinfeld case (compare \cite{H2} Example 2.7).
In order to get smooth representations with values in $\Ql$, he observes that one rather has to consider the continuous  $\ell$-adic cohomology with compact support  of these objects. The latter cohomology theory was implicitly introduced by him for adic spaces \cite{H2} (except that he does not tensor the
cohomology with $\Ql$).
In the case of schemes this construction goes back to Jannsen \cite{J}. 

We briefly recall the notion of a period domain for reductive groups \cite{RZ}. We refer to the next section for a more detailed description.
Let $L$ be an algebraically closed field of characteristic $p>0$. Denote  by $K_0=W(L)_{\mathbb Q}$ the fraction field of the ring
of Witt vectors $W(L)$ and by $\sigma$ the Frobenius automorphism of $K_0.$ 
We consider a triple
$$(G,b,\{\mu\})$$
consisting of a quasi-split reductive group $G$ defined over $\Qp,$ an element $b\in G(K_0)$ and a conjugacy class $\{\mu\}$
of 1-PS of $G.$ According to Kottwitz \cite{K1}, the element $b$ gives rise to an isocrystal with $G$-structure $N_b$  on $L.$ By definition this is an  
exact faithful tensor functor
$$N_b: Rep_{\mathbb Q_p}(G) \longrightarrow Isoc(L)$$ from the category of
finite-dimensional rational $G$-representations into the category of
isocrystals over $L.$ Here an isocrystal over $L$ is a pair $(N,\Phi)$ consisting of a finite-dimensional $K_0$-vector space $N$ and
a $\sigma$-linear bijective endomorphism $\Phi$ of $N.$ To the conjugacy class $\{\mu\}$ is associated a flag variety $\Fl=\Fl(G,\{\mu\})$  defined over some finite extension $E$ of $\Qp.$ Let $\Fl^{rig}$ be the rigid-analytic variety attached to $\Fl.$ By choosing a  $G$-representation
$V$, we get for each  point $x\in \Fl^{rig}$ a filtered isocrystal $(N_b(V),\F_x)$ over some field extension  $K$ of $K_0.$
A filtered isocrystal $(N,\Phi,\F)$ is called weakly admissible (in the sense of Fontaine \cite{Fo}) if 
$$ \sum\nolimits_i i\cdot \dim gr^i_{{\cal F}'}(N' \otimes_{K_0} K) \leq \mbox{ord}_p
\det(\Phi')$$
for all subisocrystals $(N',\Phi')$ of $(N,\Phi)$
and equality for $(N',\Phi')$ = $(N,\Phi).$ Here ${\cal F}'$ denotes
the filtration on $N'\otimes_{K_0} K$ induced by ${\cal F}.$
In \cite{RZ} it is shown that the points $x\in \Fl^{rig}$ which yield a weakly admissible
filtered isocrystal  $(N_b(V),\F_x)$ for every representation $V$ form an admissible open subset $\Fl^{wa}_b$  of  $\Fl^{rig}.$
It has therefore the structure of a rigid-analytic variety over $K_0.$ It is defined over the composite $E_s=E.\mathbb Q_{p^s},$
where  $\Q_{p^s}$ is the field of definition of $N_b$ (the unramified extension of degree $s$ of $\Qp$  in $K_0).$
The space $\Fl^{wa}_b$ is called the period domain with respect to $(G,b,\{\mu\})$. 
In this paper we consider the case of a basic isocrystal. This means that the automorphism group
$J$ of $N_b$ is an inner form of $G$. In general it is an inner form of some Levi subgroup of $G.$ 

Of general interest is the computation of the \'etale cohomology of period domains. 
The cohomology groups are equipped in a natural way with actions
of the $p$-adic Lie group $J(\Qp)$  and the absolute Galois group
$\Gamma_{E_s}=Gal(\overline{E}_s/E_s).$  In this paper we compute the continuous $\ell$-adic cohomology with compact support of $\Fl_b^{wa}$, where $\ell \in \N$ is a prime number with $(\ell,p)=1.$ The computation proceeds as follows. Consider the adic spaces $(\Fl_b^{wa})^{ad}$,
$\Fl^{ad}$ associated to $\Fl^{wa}_b$ and $\Fl.$ Let $Y^{ad}$ be the closed complement of $(\Fl_b^{wa})^{ad}$ in $\Fl^{ad}.$ In \cite{O2} we constructed an acyclic resolution of the constant \'etale sheaf $\Z/n\Z$  on  $Y^{ad}$ in terms of the Tits building of $J.$ The resulting spectral sequence degenerates  in $E_2$ and computes the
cohomology of $Y^{ad}$ with coefficients in $\Zn.$ In this paper we take up this idea and generalize this construction to
an acyclic resolution of a certain object in the derived category of \'etale sheaves on $Y^{ad}.$ 
The hypercohomology of this object appears in a long exact cohomology sequence next to the continuous $\ell$-adic cohomology groups of $\Fl_b^{wa}$ and of $\Fl.$ As in loc.cit.\! we get a spectral sequence which computes the hypercohomology of this object and hence the cohomology of $\Fl^{wa}_b.$ The cohomology groups are computed as objects in the category of smooth $J(\Qp)\times \Gamma_{E_s}$-representations.
Our result is the following:
$$H^{\ast}_c(\Fl_b^{wa},\Ql)=\bigoplus_{[w] \in W^\mu/\Gamma_{E_s}}
v^J_{P_{I_{[w]}}}\otimes
ind_{[w]}\{-l([w])\}[-|\Delta\setminus I_{[w]}|].$$
Here $W^\mu$ denotes the set of Kostant representatives in the Weyl group $W$ of $G$ with respect to the stabilizer of $\mu.$ The index set $W^\mu/\Gamma_{E_s}$ consists of the orbits for the natural action of $\Gamma_{E_s}$ on $W^\mu.$  The 
objects $v^J_{P_{I_{[w]}}}$ denote generalized Steinberg representations with respect to certain
standard parabolic subgroups $P_{I_{[w]}}$ of $J.$ On the other hand, the Galois group $\Gamma_{E_s}$
operates via Galois representations $ind_{[w]}.$ 
We refer to the next section for a detailed  explanation of this formula.

We now mention previous results on the cohomology of period domains. 
The bulk of this work has been done for Drinfeld's upper half plane $\Omega^{d+1}.$ 
The de Rham and the \'etale cohomo\-logy  with torsion coefficients prime to $p$ of $\Omega^{d+1}$  were computed by Schneider and Stuhler in \cite{SS}.
The case of the continuous $\ell$-adic cohomology with compact support has been recently solved by Dat \cite{D}.
He makes heavy use of the specialization map of $\Omega^{d+1}$ into the Bruhat-Tits building of $GL_n.$ 
For further results on the cohomology of $\Omega^{d+1}$ dealing with $p$-adic aspects, we refer to the work of Alon and de Shalit \cite{AS1}, \cite{AS2}, Gro\ss e-Kl\"onne \cite{G}, Iovita and Spiess \cite{IS},  Ito \cite{I}, Schneider \cite{S}, de Shalit \cite{dS1}, \cite{dS2}.  
\newline 
For general period domains where the underlying isocrystal is basic, there exists a formula for the continuous $\ell$-adic Euler-Poincar\'e characteristic in the Grothendieck group of smooth $J(\Qp)\times
Gal(\overline{E}_s/E_s)$-representations due to Kottwitz and Rapoport \cite{R1}, \cite{R2}. Their approach is through the Harder-Narasimhan stratification on $\Fl^{rig}$ in which the period domain appears as an open stratum. The individual cohomology groups with torsion coefficients prime to $p$ 
were computed in \cite{O2}. 

The present paper structures as follows.
In the first section we fix notation and formulate the main result. The second section deals with the construction of the 
acyclic resolution of the above mentioned object in the derived category of \'etale sheaves on $Y^{ad}.$ In the final section
we compute the cohomology by evaluating the spectral sequence associated to this complex. 

I wish to thank M. Rapoport for his useful comments on this paper.
Furthermore, I would like to thank J.-F. Dat and A. Huber for helpful and interesting discussions.

\section{Notation and the main result}
The notation will be the same as in \cite{O2}. So, $L$ is an algebraically closed  field of characteristic $p>0,$ 
$K_0=W(L)_{\mathbb Q}$ denotes the fraction field of the ring
of Witt vectors and  $\sigma \in  \mbox{Aut}(K_0/\mathbb Q_p)$ is the Frobenius
automorphism.  Furthermore, $\ell\in \N$ is a prime number with $(\ell,p)=1.$
 
Let $G$ be a connected quasi-split reductive group over $\Qp.$
We fix an element $b\in G(K_0)$ which defines an isocrystal with
$G$-structure on $L$, i.e, an exact faithful tensor functor
$$N_b: Rep_{\mathbb Q_p}(G) \longrightarrow Isoc(L)$$ from the category of
finite-dimensional rational $G$-representations into the category of
isocrystals over $L.$  We denote by
$$\nu:=\nu_b: \mathbb D \longrightarrow G_{K_0}$$ 
the slope homomorphism of $N_b,$ where $\mathbb D$ is the algebraic pro-torus over $\Qp$ with character
group $\mathbb Q.$ For a representation $\rho : G \rightarrow GL(V),$ the composition
$\rho\circ \nu:\mathbb D \rightarrow GL(V_{K_0})$ describes the decomposition of $N_b(V)$ into simple subisocrystals together with their slopes (\cite{RR} Example 1.10).
We suppose that $b$ is  decent and basic. 
By definition (\cite{RZ} 1.8), a decent  element $b\in G(K_0)$ satisfies for some integer $s\in \N$ an equation
$$(b\sigma)^s=s\nu_b(p)\sigma^s $$
in the semi-direct product $G(K_0) \rtimes \langle \sigma \rangle$
 such that $s\cdot\nu_b:\mathbb D
\longrightarrow G_{K_0}$ factors through $\mathbb G_m.$
It follows that $b\in G({\mathbb Q}_{p^s})$ and that $\nu$ is defined over ${\mathbb Q}_{p^s}.$
The element $b$ is basic if $\nu$ factors through the center of $G_{K_0}.$ 
An alternative definition is as follows.
Let $J$ be the automorphism group of $N_b,$ which is a reductive algebraic group over $\Qp$. 
In general it is an inner form of some Levi subgroup of $G.$ It is
is an inner form of G if and only if $b$ is basic.

We fix a conjugacy class $$\{\mu\} \subset X_\ast(G)$$ of one-parameter
subgroups (1-PS) of $G$ over the algebraic closure $\overline{\mathbb Q}_p$  of $\Qp.$ 
Let $E$ be the Shimura field of $\{\mu\}.$ This 
finite extension of $\Qp$ is the smallest extension of $\Qp$ such that its absolute Galois group
stabilizes $\{\mu\}$ as an element in the set of conjugacy classes of 1-PS of $G.$ By Kottwitz's Lemma
(\cite{K2} Lemma 1.1.3), there exists a 1-PS $\mu\in \{\mu\}$ which is defined over $E.$
Hence, the conjugacy class $\{\mu\}$ defines  a flag variety 
$$\Fl:= \Fl(G,\{\mu\}):=G_E/P(\mu)$$ 
defined over $E.$
Here, we write $P(\lambda)$ for 
the parabolic subgroup of $G$ which is attached to a 1-PS $\lambda.$
Let $\Fl^{rig}$ be the rigid analytic flag variety associated to $\Fl.$ Recall that its
underlying set is given by $\Fl(\Cp),$
where $\C_p:=\hat{\overline{\mathbb Q}}_p$ is the $p$-adic completion of 
$\overline{\mathbb Q}_p.$
Put $$E_s:=E.{\mathbb Q}_{p^s} \mbox{ and } \Gamma_{E_s}:=Gal(\overline{E}_s/E_s).$$
The period domain $\Fl^{wa}_b$ associated to the triple $(G,b,\{\mu\})$ is the set of filtrations in $\Fl^{rig}$ which are
weakly admissible with respect to $N_b.$
It has a natural structure of an admissible open rigid-analytic subset of $(\Fl
\otimes _E E_s)^{rig}$ (\cite{RZ} Proposition 1.36). Furthermore, it is a $J(\Qp)$-invariant subset of $\Fl^{rig}$ with respect 
to the natural action of $J(\mathbb Q_p)$ on $\Fl^{rig}.$

We fix a maximal $\Qp$-split torus $S$ of the derived group
$J_{der}$ of $J.$ Further, we let $P_0$ be a minimal $\Qp$-parabolic subgroup of $J$  such
that $S\subset P_0.$ We denote by
$$\Delta=\{\alpha_1,\ldots,\alpha_d \} \subset X^\ast(S)$$ 
the corresponding set of relative simple roots.
Let $$\{\omega_\alpha; \alpha \in \Delta \}\subset
X_\ast(S)_{\mathbb Q}$$ be the dual basis of $\Delta$ 
with respect to the natural pairing $X_\ast(S)\times X^\ast(S) \rightarrow \Z.$
Let $T$ be a maximal torus  of $G$ which contains $S$ and such that $\mu,\nu \in
X_\ast(T)_{\Q}\cong \mbox{Hom}_{K_0}(\mathbb D,T)$. 
We fix a Borel subgroup $B$ of $G$ such that
\begin{itemize}
\item $B\subset P(\omega_\alpha)$ for all $\alpha \in \Delta$ 
\item $\mu$ lies in the positive Weyl chamber with respect to $B.$
\end{itemize}
We may assume that $\Delta$ is given by restriction of a  root basis of $G$ with respect to $B \supset T.$ 
Put 
$$\bar{\mu}:=\frac{1}{|Gal(E/\Qp)|}\sum_{\gamma\in Gal(E/\Qp)}\gamma \mu.$$
By \cite{FR} we know that $\Fl_b^{wa}$ is non-empty if and only if 
$$\bar{\mu} \geq \nu$$ 
with respect  to the dominance order $\geq$ on $X_\ast(T)_{\mathbb Q}$  induced by $B.$ In our situation of a basic isocrystal, the latter condition simply means that $\bar{\mu}- \nu \in X_\ast(T_{der}),$
where $T_{der}\subset T$ is the maximal torus of the derived group $G_{der}$ of $G.$ 

Let $W$ the Weyl group of $G.$ 
Denote by $W_\mu \subset W$ the
the stabilizer of $\mu,$ and by  $W^\mu$ the set of Kostant representatives with respect
to $W/W_\mu.$ We have an action of $\Gamma_{E_s}$ on $W$ 
which preserves $W^\mu.$ 
For any orbit $[w]\in W^\mu/\Gamma_{E_s},$ we have the induced representation 
$$ind_{[w]}:=Ind_{Stab_{\Gamma_{E_s}}(w)}^{\:\Gamma_{E_s}}(\Ql)$$
of $\Gamma_{E_s},$ where we consider the trivial action of $Stab_{\Gamma_{E_s}}(w)$ on $\Ql.$
For any subset $I \subset \Delta,$ we set
$$\Omega_I:=\Big\{[w] \in W^\mu/\Gamma_{E_s};\; (w\mu,\omega_\alpha) > (\nu,\omega_\alpha)  \; \forall \alpha \not\in I\Big\}.$$
Here, we have fixed an invariant inner product $(\; ,\;)$ on $G$ (see \cite{T} section 7). 
Finally, we denote for $[w] \in W^\mu/\Gamma_{E_s},$  by $I_{[w]}$ the
minimal subset of $\Delta$ such that $[w]$ is contained in
$\Omega_{I_{[w]}}.$  

For a parabolic subgroup $P\subset J$ defined over $\Qp,$ let 
$$i^J_P=C^\infty(J(\Qp)/P(\Qp),\Ql)$$
be the smooth representation of $J(\Qp)$ consisting of locally constant
$\Ql$-valued functions on the $p$-adic manifold $J(\Qp)/P(\Qp).$
The generalized Steinberg representation  $v^J_{P}$ is defined by 
$$v^J_P=i^J_P/\sum_{P\subsetneq Q}i^J_Q.$$  
This is an irreducible smooth admissible representation of $J(\Qp)$ (\cite{C} Theorem 1.1).
Finally, for $I\subset \Delta$ we
put 
$$P_I:=\bigcap_{\alpha \not\in I }P^J(\omega_\alpha)$$ which is a standard-parabolic subgroup of $J$
defined over $\Qp.$ The extreme cases are $P_\Delta=J$ and $P_\emptyset =P_0.$
Here, $P^J(\lambda)$ is the parabolic subgroup of $J$ attached to a 1-PS $\lambda$ of $J.$

We denote by $H^\ast_c(\Fl_b^{wa},\Ql)$  the continuous $\ell$-adic cohomology with compact support of $\Fl_b^{wa}$. We refer to section 3 for its definition.
We have a natural action of  $J(\Qp) \times \Gamma_{E_s}$ on $H^\ast_c(\Fl_b^{wa},\Ql)$ which is smooth in the sense of $p$-adic Lie group representations. In fact, in the case of finite cohomology coefficients this statement is due to Berkovich \cite{B1}. In our situation, there are unpublished notes \cite{B2} of Berkovich
treating our case. In any case, we first follow \cite{B1} Prop. 6.4 to deduce that the action of $J(\Qp)\times \Gamma_{E_s}$ on $\Fl_b^{wa}$ is continuous. Then we apply Fargue's result \cite{F} Cor. 4.1.19 to conclude
that $H^\ast_c(\Fl_b^{wa},\Ql)$ is a smooth representation.

Our result is the following theorem.
\begin{Theorem} Let $b\in G(K_0)$ be a basic element such that $\bar{\mu} -\nu \in X_\ast(T_{der}).$  We have
for the compactly supported continuous $\ell$-adic cohomology  of $\Fl_b^{wa}$ the following formula as representations of $J(\Qp) \times \Gamma_{E_s}$:
$$H^{\ast}_c(\Fl_b^{wa},\Ql)=\bigoplus_{[w] \in W^\mu/\Gamma_{E_s}}
v^J_{P_{I_{[w]}}}\otimes
ind_{[w]}\{-l([w])\}[-|\Delta\setminus I_{[w]}|] .$$
\end{Theorem}
Here for a representation $V,$ we have set $V\{-m\}=V(-m)[-2m], m\in \mathbb Z .$ The symbol $(m)$ denotes the $m$-th Tate twist and $[-m]$ symbolizes that the
corresponding module is shifted into degree $m$ of the graded cohomology ring.

\section{The fundamental complex}

In \cite{O2} we constructed
an acyclic complex of \'etale sheaves on the complement of $\Fl_b^{wa}$ in $\Fl.$ Strictly speaking,
we formed the complement in the category of Huber's adic spaces \cite{H1}, i.e., we  set
$$Y^{ad}:=\Fl^{ad}\setminus (\Fl_b^{wa})^{ad}.$$ 
Here, $^{ad}$ indicates the adic space attached to a scheme or a rigid-analytic variety.
The advantage of this larger category is that $Y^{ad}$
has the structure of a closed pseudo-adic subspace  of $\Fl^{ad}$ for which 
there exists an \'etale site and hence a topos as well (compare \cite{H1} chapter 2.3).
We write $X_{\acute{e}t}$ for the \'etale site of a scheme, a rigid-analytic variety or that of a pseudo-adic space $X.$

In \cite{H2} (see chapter 1) Huber introduced the compactly supported cohomology of $\ell$-adic sheaves on adic spaces. 
We briefly explain Huber's  definition.
Let $R$ be a complete discrete valuation ring with maximal ideal ${\frak m}$ and with char$(R/\frak m)>0.$
Let $X$ be a taut separated pseudo-adic space locally of $^+$weakly finite type over an analytic algebraically closed affinoid field $Spa(k,k^0)$. 
A $R_\bullet-$module on $X_{\acute{e}t}$ is a projective system
$$\cdots \rightarrow F_{n+1} \rightarrow F_n \rightarrow \cdots \rightarrow F_1 $$
of $R$-modules on $X_{\acute{e}t}$ such that ${\frak m}^n\cdot F_n=0$ for every $n\in \mathbb N.$ 
Denote by $mod(X_{\acute{e}t}-R_\bullet)$
the category of $R_\bullet$-modules on $X_{\acute{e}t}.$ This is an abelian category with enough injectives.
Further, we let $mod(R)$ resp. $mod(X_{\acute{e}t}-R)$ be the category of $R$-modules resp. $R$-modules on $X_{\acute{e}t}.$  
For any $R_\bullet-$module 
$(F_n)_n$ on $X_{\acute{e}t}$, the compactly supported cohomology $ H^p_c(X,(F_n)_n)$ is defined 
as follows. In the case where
$X$ is partially complete, the cohomology is given by
$$H^p_c(X,(F_n)_n):=R^p\Gamma_c(X,(F_n)_n), \;p\in \N.$$
Here $\Gamma_c$ is the left exact functor
\begin{eqnarray*} 
\Gamma_c:mod(X_{\acute{e}t}-R_\bullet) & \rightarrow & mod(R) \\
(F_n)_n & \mapsto & \Gamma_c(X,\varprojlim F_n).
\end{eqnarray*}
where $\Gamma_c(X,\varprojlim F_n)$ is the $R$-module of global sections whose support is complete
over $Spa(k,k^0).$ In general, if $X$ is not partially complete, one defines the compactly supported cohomology  using a partially complete compactification of $X.$ 
Finally, by the continuous $\ell$-adic cohomology with compact support of $X$ we mean the $\Ql$-vector spaces
$$H^\ast_c(X,\Ql)= H^\ast_c(X, (\Zln)_n)\otimes_{\Zl} \Ql.$$
In the case of algebraic varieties this definition is due to Jannsen \cite{J}. 
Huber does not use this notation in \cite{H2}.

We indicate briefly the relation to Berkovich's definition of compactly supported cohomology of $R_\bullet$-modules which applies to his analytic spaces. Berkovich defines this cohomology  by using a different global section functor. It is given by 
\begin{eqnarray*}
mod(X_{\acute{e}t}-R_\bullet) &\rightarrow& mod(R) \\ \\
(F_n)_n &\mapsto& \varinjlim_{U\in \mathbb U} \varprojlim_{n\in \N} \Gamma_c(U,F_n).
\end{eqnarray*}
Here $\mathbb U$ is the set of distinguished open subsets in the analytic space $X$.
In his paper (\cite{F} 4.1.2), Fargues uses Huber's definition of compactly supported cohomology of $R_\bullet$-modules in the case of analytic spaces. He proves that it coincides with that of Berkovich. For rigid-analytic varieties, this statement was already shown in \cite{H2} Prop.\! 1.5.

Consider the left exact functors 
\begin{eqnarray*}
\pi_\ast:mod(X_{\acute{e}t}-R_\bullet) & \rightarrow & mod(X_{\acute{e}t}-R) \\
(F_n)_n & \mapsto & \varprojlim_n F_n 
\end{eqnarray*}
and
\begin{eqnarray*}
\Gamma_!:mod(X_{\acute{e}t}-R) & \rightarrow & mod(R)  \\
F & \mapsto & \Gamma_c(F).
\end{eqnarray*}
If $X$ is partially complete, then,  by definition, we have $\Gamma_c=\Gamma_! \circ \pi_\ast.$
This identity then gives rise to an identity (compare \cite{H2} Lemma 2.3 (i))
\begin{equation}
R^+\Gamma_c=R^+\Gamma_! \circ R^+\pi_\ast
\end{equation}
of exact functors from the derived category $D^{\geq 0}$(mod$(X_{\acute{e}t}-R_\bullet))$ into  $D^{\geq 0}($mod$(R)).$

\bigskip
For any rational $1$-PS $\lambda \in X_\ast(J)_{\Q}$ of $J,$ we consider in \cite{O2}
the closed $E_s$-subvariety $Y_\lambda$ of $\Fl$ consisting of points where $\lambda$
violates the semi-stability condition.
By \cite{O2} Corollary 2.4 we have the identity
\begin{eqnarray} 
Y^{ad}= \bigcup_{\lambda}Y_\lambda^{ad}.
\end{eqnarray}
of pseudo-adic spaces.
Furthermore, we define for any subset $I\subset \Delta$ 
the closed subvariety 
$$ Y_I:=\bigcap_{\alpha \notin I} Y_{\omega_\alpha}$$
of $\Fl,$ on which we have an action of the $p$-adic group $P_I(\Qp).$ 
For any compact open subset $ W\subset J/P_I(\Qp),$ we put
$$Z_I^W:= \bigcup_{g\in W} gY_I^{ad}.$$
Then by \cite{O2} Lemma 3.2 we know 
that $Z_I^W$ is a closed pseudo-adic subspace of $\Fl^{ad}.$
By (2) it follows that
$$Y^{ad}=\bigcup_{I\subset \Delta \atop |\Delta\setminus I|=1}Z_I^{J/P_I(\Qp)}.$$

Starting from a subset $I\subset \Delta$ and the constant \'etale sheaf $\Zn ,\;n\in \N,$ on $Y^{ad}$ we  constructed in loc.cit.\! a sheaf of locally constant sections on the same space. Actually, the
same construction goes through if we start with an arbitrary \'etale sheaf $F$ on $Y^{ad}.$ We recall  the definition. 
Let
$$\Phi_{g,I} :gY_I^{ad} \longrightarrow Y^{ad}$$
resp.
$$\tilde{\Phi}_{g,I,W} :gY_I^{ad} \longrightarrow Z_I^W$$
resp.
$$\Psi_{I,W} :Z^W_I \longrightarrow Y^{ad}$$
the natural closed embeddings of pseudo-adic spaces. 
Let $F$ be a $R$-module on $Y^{ad}_{\acute{e}t}.$
Put 
$$ F_{g,I}:=(\Phi_{g,I})_\ast(\Phi_{g,I}^\ast F)$$ 
resp. 
$$ F_{Z_I^W}:=(\Psi_{I,W})_\ast(\Psi_{I,W}^\ast F)$$
and let
$$\tilde{\Phi}_{g,I,W}^{\#} :F_{Z_I^W} \longrightarrow F_{g,I}$$
be the natural homomorphism given by restriction.
We denote by
$$\sideset{}{'}\prod_{g\in J/P_I (\Qp)}F_{g,I} $$ the $R$-submodule of 
$\prod_{g\in J/P_I(\Qp)} F_{g,I},$ which is defined as the sheaf associated to the following presheaf $\;\;^P\prod\nolimits_{g\in J/P_I (\Qp)}' F_{g,I}.$
For any element $U \rightarrow Y^{ad} $ of the \'etale site $Y^{ad}_{\acute{e}t},$ we put

\noindent $\begin{array}{lcl} \Big(\displaystyle\; ^P\!\!\!\!\!\sideset{}{'}\prod_{g\in J/P_I(\Qp)} F_{g,I}\Big)(U) \!\!\! & := & \!\! \Big\{ (s_g)_g \in \displaystyle \prod_{g\in J/P_I(\Qp)} 
F_{g,I} (U);   \mbox{ there exists a (finite) disjoint cover- }  \\
& & \mbox{ing } J/P_I(\Qp)= \bigcup\limits_{j\in A}^\cdot W_j  \mbox{ by compact open  subsets  and sections }\\ 
& & \mbox{} s_j \in F_{Z_I^{W_j}}(U) ,\; j \in A, \mbox{such that } \tilde{\Phi}^{\#}_{g,I,W_j}(s_j)= s_g \; \mbox{for all } g\in W_j \Big\}.  \\ 
\end{array}$

\noindent We call $\prod\nolimits'_{g\in J/P_I(\Qp)} F_{g,I}$ the subsheaf of locally constant sections of $\prod\nolimits_{g\in J/P_I(\Qp)} F_{g,I}$.

\noindent {\bf Remark:} If we work with the restricted \'etale site
$Y^{ad}_{\acute{e}t,f.p.}$ consisting of objects $U$ in $Y^{ad}_{\acute{e}t}$ where
the structure morphism $U\rightarrow Y^{ad}$ is quasi-compact and quasi-separated, it is easy to see that the
presheaf $\;\;^P\prod\nolimits_{g\in J/P_I (\Qp)}' F_{g,I}$ is already a sheaf (compare also \cite{O2}). 

The sheaf $\prod\nolimits'_{g\in J/P_I(\Qp)} F_{g,I}$ may be viewed as an
inductive limit of sheaves. In fact, let 
${\cal C}_I$ be the category of compact open disjoint coverings of $J/P_I(\Qp)$
in which the morphisms are given by the refinement-order.
For a covering $c =(W_j)_{j\in A} \in {\cal C}_I,$ we denote by $F_c$ the sheaf defined by

\noindent $\begin{array}{lcl}F_c(U) & := & \Big\{ (s_g)_g \in \displaystyle\prod_{g\in J/P_I(\Qp)} F_{g,I}(U) ; 
\;\mbox{ there are sections }\; s_j\; \in F_{Z_I^{W_j}}(U), j\in A, \mbox{ such } \\
& &  \mbox{ that } \tilde{\Phi}^{\#}_{g,I,W_j}(s_j) = s_g \; \mbox{for all } g \in W_j \Big\} .
\end{array} $

\noindent Alternatively, one can define
$F_c$ as the image of the natural morphism of sheaves
\begin{equation}
\bigoplus_{j \in A}\; F_{Z_I^{W_j}} \hookrightarrow \prod_{g\in J/P_I(\Qp)} F_{g,I}.
\end{equation}
Then we may write 
\begin{eqnarray}
\sideset{}{'}\prod_{g\in J/P_I(\Qp)} F_{g,I} = \varinjlim_{c\in {\cal C}_I} F_c.
\end{eqnarray}

Our construction extends naturally to the derived category of \'etale sheaves on $Y^{ad}.$
Indeed, the assignment $F\mapsto \prod\nolimits'_{g\in J/P_I(\Qp)} F_{g,I}$ is functorial and induces therefore a functor
$$\displaystyle \sideset{}{'} \prod_{g\in J/P_I(\Qp)}\!\!\!\!\!\!\! : mod(Y^{ad}_{\acute{e}t}-R) \rightarrow mod(Y^{ad}_{\acute{e}t}-R),$$
which is easily seen to be exact. Thus, we obtain an exact functor on the corresponding derived categories
$$\displaystyle \sideset{}{'}\prod_{g\in J/P_I(\Qp)} \!\!\!\!\!\!\!= R^+ \sideset{}{'}\prod_{g\in J/P_I(\Qp)} : D^{\geq 0}(mod(Y^{ad}_{\acute{e}t}-R)) \rightarrow D^{\geq 0}(mod(Y^{ad}_{\acute{e}t}-R)).$$

\bigskip
From now on we let $R=\Zl.$
Furthermore, we consider the ${\Zl}_\bullet$-module 
$$(F_n)_n=(\Zln)_n \in mod(\Fl^{ad}-{\Zl}_\bullet).$$
Let $i:Y^{ad} \hookrightarrow \Fl^{ad}$ be the inclusion. Put
$$ \Upsilon :=i^\ast R^+\pi_\ast(\Zln)_n \in D^{\geq 0}(mod(Y^{ad}_{\acute{e}t}-\Zl)).$$
By \cite{H2} Prop. 2.6 (ii), we have a long exact sequence
\begin{eqnarray}
\cdots  &\rightarrow &  H^p_c((\Fl^{wa})^{ad},(\Zln)_n) \rightarrow H^p_c(\Fl^{ad},(\Zln)_n) \rightarrow
 H^p_c(Y^{ad},\Upsilon) \nonumber \\ \nonumber \\
&\rightarrow &   H^{p+1} _c((\Fl^{wa})^{ad},(\Zln)_n) \rightarrow \cdots
\end{eqnarray}
Taking Huber's comparison theorem for schemes into account (\cite{H2} Theorem 4.2) 
we have 
$$H_c^\ast(\Fl^{ad},(\Zln)_n) = H^\ast_c(\Fl,(\Zln)_n).$$ 
On the other hand, by Jannsen's comparison theorem (\cite{J} Remark 3.5) we conclude that
$$H^\ast_c(\Fl,(\Zln)_n) = H^\ast_c(\Fl,\Zl)=H^\ast(\Fl,\Zl),$$
where   $H^\ast(\Fl,\Zl):=\varprojlim_n H^\ast(\Fl,\Zln)$ denotes the $\ell$-adic cohomology of
the scheme  $\Fl.$
Further, by the very definition of the compactly supported $\ell$-adic cohomology of rigid-analytic varieties we have
$$ H^\ast_c((\Fl^{wa})^{ad},(\Zln)_n)=H^\ast_c(\Fl^{wa},(\Zln)_n).$$
Therefore, in order to prove Theorem 1.1 it suffices to compute the cohomology groups
$$H^p_c(Y^{ad},\Upsilon)=H^p(Y^{ad},\Upsilon).$$

Let $F$ be a \'etale sheaf on $Y^{ad}$. In \cite{O2} we constructed in the case $F=\Zn$ a complex of \'etale sheaves on $Y^{ad}:$

\begin{equation}\tag{$*$}
\begin{aligned} 0 \rightarrow F \rightarrow\!\!\! \bigoplus_{I \subset
\Delta \atop |\Delta\setminus I|=1} \sideset{}{'}\prod_{g\in J/P_I(\Qp)} F_{g,I}
\rightarrow \!\!\!\bigoplus_{I \subset \Delta \atop |\Delta\setminus I|=2}
\sideset{}{'}\prod_{g\in J/P_I(\Qp)} F_{g,I} \rightarrow \cdots  
\\ \\
 \dots \rightarrow
\!\!\!\bigoplus_{I \subset \Delta \atop |\Delta\setminus I|=d-1} \sideset{}{'}\prod_{g\in J/P_I(\Qp)} F_{g,I}
\rightarrow \sideset{}{'}\prod_{g\in J/P_\emptyset(\Qp)} F_{g,\emptyset} \rightarrow 0
\end{aligned}
\end{equation}

\noindent The same construction works for  arbitrary \'etale sheaves $F$ on $Y^{ad}$ as well.  Furthermore,  we  proved 
(loc.cit. Theorem 3.3) that
the complex above is acyclic in the case $F=\Z/n\Z$. The same holds if we consider overconvergent 
\'etale sheaves $F$ on $Y^{ad}.$ Recall that an \'etale sheaf $F$ on a pseudo-adic space $X$ is called overconvergent
if the specialization morphisms $F_{\xi_1} \rightarrow F_{\xi_2}$ are isomorphisms for all pairs of geometric points ${\xi_1},\xi_2$ of $X$ such that $\xi_1$ is a specialization of $\xi_2.$

\begin{Theorem} Let $F$ be an overconvergent \'etale sheaf on $Y^{ad}$. Then the complex $(*)$ is acyclic.
\end{Theorem}

\proof The proof is almost the same as that of Theorem 3.3 in \cite{O2}. First, we follow the proof
of Lemma 3.4 in loc.cit.\! to show that the complex $(*)$ consists of overconvergent sheaves.
Further, by localizing in a maximal geometric point $\xi$ of $Y^{ad}$ we obtain a chain complex with values in $F_\xi$ instead of $\Zn$. In loc.cit. we reduced the acyclicity of this chain complex in
the case $F=\Z/n\Z$ to the contractibility of its index set. The same argument applies in our general case. \qed

\bigskip
\noindent Choose a complex $\overline{\Upsilon}$ representing the object $\Upsilon \in D^{\geq 0}(mod(Y^{ad}_{\acute{e}t}-\Zl)$ .  Substituting the \'etale sheaf $F$ in $(*)$ by $\overline{\Upsilon}$ we obtain
a  complex in the category of complexes with values in $mod(Y^{ad}_{\acute{e}t}-\Zl)$:

\begin{equation} \tag{$**$} 
\begin{aligned}
0 \rightarrow \overline{\Upsilon}  \rightarrow\!\!\! \bigoplus_{I \subset
\Delta \atop |\Delta\setminus I|=1} \sideset{}{'}\prod_{g\in J/P_I(\Qp)} \overline{\Upsilon}_{g,I}
\rightarrow \!\!\!\bigoplus_{I \subset \Delta \atop |\Delta\setminus I|=2}
\sideset{}{'}\prod_{g\in J/P_I(\Qp)} \overline{\Upsilon}_{g,I} \rightarrow \dots 
\\ \\
\cdots \rightarrow
\!\!\!\bigoplus_{I \subset \Delta \atop |\Delta\setminus I|=d-1} \sideset{}{'}\prod_{g\in J/P_I(\Qp)} \overline{\Upsilon}_{g,I} \rightarrow \sideset{}{'}\prod_{g\in J/P_\emptyset(\Qp)} \overline{\Upsilon}_{g,\emptyset} \rightarrow 0.
\end{aligned}
\end{equation}

\begin{Corollary}
The complex $(**)$ is exact.
\end{Corollary} 

\proof By Theorem 2.1 we merely have to assure that the individual contributions of $\Upsilon=i^\ast R^+\pi_\ast(\Zln)_n$ are overconvergent. Since $Y^{ad}$ is a closed subspace of $\Fl^{ad}$, it is enough to prove the assertion for the complex $R^+\pi_\ast(\Zln)_n$ on $\Fl^{ad}.$
The \'etale topoi of $\Fl^{ad}$ and $\Fl^{rig}$ are equivalent \cite{H1} Prop.\! 2.1.4. 
Let $\Fl^{an}$ be Berkovich's analytic space attached to $\Fl.$
There is an equivalence between the \'etale topos of $\Fl^{an}$ and the \'etale subtopos of $\Fl^{rig}$ consisting of overconvergent sheaves (loc.cit. Prop.\! 8.3.3). 
The claim follows now easily by comparison.
\qed

\section{The computation}
In this final section we evaluate the spectral sequence which is induced by
the complex $(**)$. This computation ends with the
determination of the hypercohomology of $\Upsilon.$ 
The proof of Theorem 1.1 then follows from the long exact cohomology
sequence (5). The arguments of this final section are modelled on the results in \cite{O2} section 4
treating the torsion sheaf $\Zn$ instead of $\Upsilon.$

For the next two statements we fix a subset $I\subsetneq \Delta.$ 

\begin{Proposition} We have 
\begin{equation} H_c^\ast(Y_I^{ad},\Ql)=\bigoplus_{[w] \in \Omega_I}
  ind_{[w]}\{-l[w]\}
\end{equation}
\end{Proposition}
\proof By the comparison theorems of Huber (see \cite{H2} Thm. 4.2)  resp. Jannsen (see \cite{J} Remark 3.5), we have $H_c^\ast(Y_I^{ad},(\Zln)_n)= H_c^\ast(Y_I,(\Zln)_n)$ resp. 
$H_c^\ast(Y_I,(\Zln)_n)=H^\ast(Y_I,\Zl).$
Then we proceed as in Proposition 4.2  \cite{O2}. \qed

\noindent Next we compute the cohomology of the ingredients
in $(**)$, i.e., the hypercohomology of the objects $\prod\nolimits'_{g\in J/P_I(\Qp)}\Upsilon_{g,I}.$ 

\begin{Proposition} We have 
$$H^i(Y^{ad}, \sideset{}{'}\prod_{g\in J/P_I(\Qp)} \Upsilon_{g,I}) =
  C^\infty(J/P_I(\Qp),  H^i (Y_I^{ad},\Zl)) $$
for all i $\in \N.$ Here,  $C^\infty(J/P_I(\Qp),  H^i (Y_I^{ad},\Zl))$ denotes the space of locally constant functions on $J/P_I(\Qp)$ with values in   $H^i(Y_I^{ad},\Zl).$
\end{Proposition} 

\proof The proof is similar to Prop.\! 4.3 in loc.cit. Again, we fix as in the construction
of $(**)$ a representative $\overline{\Upsilon}$ of $\Upsilon.$ First of all, we have by (4) the identity
$$\sideset{}{'}\prod_{g\in J/P_I(\Qp)} \overline{\Upsilon}_{g,I}
=\varinjlim_{c\in {\cal C}_I} \overline{\Upsilon}_c $$ 
Since $Y^{ad}$ is quasi-compact and  ${\cal C}_I$ is 
pseudo-filtered, we get (compare also \cite{H1} 2.3.13)
$$ H^i(Y^{ad}, \sideset{}{'}\prod\limits_{g\in J/P_I(\Qp)} \Upsilon_{g,I}) =
\varinjlim_{c\in {\cal C}_I} H^i(Y^{ad},\Upsilon_c).$$ 
But by (2) we have a decomposition  $\overline{\Upsilon}_c \cong \bigoplus_{W\in c} \overline{\Upsilon}_{Z^W_I}.$ Thus, we obtain
$$H^i(Y^{ad}, \Upsilon_c)= \bigoplus_{W\in c} H^i(Z^W_I,\Upsilon_{|Z^W_I}). $$ 
Let $(W_s)_{s\in \N}$ be a cofinal family of compact open neighbourhoods of
the base point $1\cdot P_I$ in the $p$-adic manifold $J/P_I(\Qp).$ 
By Lemma 4.4 of loc.cit. we have the identity
$$\bigcap_{s \in \N} {Z_I^{W^s}} = Y_I^{ad}.$$ 
Now, we apply \cite{H1} 2.4.6\footnote{Actually the statement there is only formulated for sheaves. But it is easy to see that the argument generalizes to the derived category} to conclude that
$$\varinjlim_{s \in \N} H^i(Z_I^{W^s},\Upsilon_{|Z^{W^s}_I} ) 
= H^i(Y^{ad}_I,\Upsilon_{|Y^{ad}_I}).$$
Using (1) we see that
$$H^i(Y^{ad}_I,\Upsilon_{|Y^{ad}_I})=H^i(Y_I^{ad},(\Zln)_n).$$
so that by reapplying Huber's comparison theorem (\cite{H2} Thm. 4.2) we obtain 
$$H^i(Y^{ad}_I,\Upsilon_{|Y^{ad}_I})=H^i(Y_I^{ad},\Zl).$$
Combining these facts we get
\begin{eqnarray*}
H^i(Y^{ad}, \sideset{}{'}\prod_{g\in J/P_I(\Qp)} \Upsilon_{g,I})  
=  \varinjlim_{c\in {\cal C}_I} \bigoplus_{W\in c} H_c^i(Z^W_I,\Upsilon_{|Z^W_I} ) =   C^\infty(J/P_I(\Qp),  H^i(Y_I^{ad},\Zl)) . 
\end{eqnarray*} \qed

Consider the spectral sequence  induced by the acyclic complex $(**):$
$$E_1^{p,q} = H_c^q(Y^{ad},\bigoplus_{I \subset \Delta
\atop |\Delta\setminus I|=p+1} \sideset{}{'}\prod_{g\in J/P_I(\Qp)}\Upsilon_{g,I})\otimes_{\Zl} \Ql \Longrightarrow
H_c^{p+q}(Y^{ad},\Upsilon)\otimes_{\Zl} \Ql.$$ 
From the previous proposition we deduce that
$$ E^{p,q}_1= \bigoplus_{I\subset \Delta \atop |\Delta\setminus I|=p+1} C^\infty(J/P_I(\Qp),  H^q(Y_I^{ad},\Ql)).$$ 
As in \cite{O2}, we obtain by Proposition 3.1
  a decomposition into subcomplexes
$$E_1 =\bigoplus_{[w]\in W^{\mu}/\Gamma_{E_s}} E_{1,[w]},$$
where $E_{1,[w]}$ is the complex
$$\Big(\bigoplus_{I_{[w]} \subset I \atop |\Delta\setminus I|=1}i^J_{P_I} \otimes ind_{[w]}
\rightarrow \bigoplus_{I_{[w]} \subset I \atop |\Delta\setminus I| =2}
i^J_{P_I} \otimes ind_{[w]} \rightarrow \cdots  
\rightarrow i^J_{P_{I_{[w]}}}\otimes ind_{[w]})\Big)\{-l([w])\}.$$
By applying Prop.\! 4.5 of loc.cit.\! , we obtain for the second stage of the spectral sequence $E_2=\bigoplus_{[w]\in W^\mu/\Gamma_{E_s}} E_{2,[w]}$ the following expression:

$$ E^{0,r}_{2,[w]}=\left\{
\begin{array}{l@{\quad  \quad}l}
   i^J_{P_{I_{[w]}}}\otimes ind_{[w]}(-l([w]))  &:\;  |\Delta\setminus I_{[w]}|=1 ,\;2l([w])=r \\ \\
i^J_J\otimes ind_{[w]}(-l([w])) &:\; |\Delta \setminus I_{[w]}|>1,\; 2l([w])=r \\ \\
0 &:\; otherwise
\end{array}
\right.  $$

\noindent respectively for $p>0$
$$ E^{p,r-p}_{2,[w]}=\left\{
\begin{array}{l@{\quad  \quad}l}
v^J_{P_{I_{[w]}}}\otimes ind_{[w]}(-l([w]))  &:\;p=|\Delta \setminus I_{[w]}|-1,\; 2l([w]) + p = r \\ \\
0 &:\; otherwise
\end{array} .
\right.  $$
\noindent Since $E_s$ is a local field (Beware of the double-notation !), we conclude by
weight arguments (compare also the corresponding argument in loc.cit. section 5) that the 
spectral sequence degenerates in $E_2.$ 
Thus, we have for $p,r\in \N,$  

\begin{equation}
gr^p(H_c^r(Y^{ad},\Upsilon)\otimes_{\Zl}\Ql)=E_\infty^{p,r-p}=E_2^{p,r-p}=\bigoplus_{[w]\in W^\mu/\Gamma_{E_s}} E_{2,[w]}^{p,r-p}.
\end{equation}


In order to show that the canonical filtration on $E_\infty$ splits, we proceed as follows. 
By unpublished results of Berkovich respectively by Dat (\cite{D} appendix B.2) we know that the cohomology groups $H_c^r(Y^{ad},\Upsilon)$ are smooth $J(\Qp)$-modules. In  \cite{O1} (resp. \cite{D} in the split case) the following statement on extensions of generalized Steinberg representations was proved.

\begin{Theorem} \cite{O1} Let $r\in \N$ be the $\Qp$-rank of the center of $J.$
Let $I,I'\subset \Delta.$ Then for any $i \in \N$ with $r\geq i \geq |I\cup
  I'| -  |I\cap I'|,$
$$Ext_{J(\Qp)}^i(v^J_{P_I},v^J_{P_{I'}})=  \Ql^{r \choose j },\; \mbox{ where }   j= i- |I\cup
  I'| -  |I\cap I'| .$$
For all other $i$ the Ext-group vanishes.
\end{Theorem}

\noindent Consider the equation  
$$2l([w]) + |\Delta\setminus I_{[w]}|-1= r= 2l([w']) +
|\Delta\setminus I_{[w']}|-1$$
\noindent with $[w],[w'] \in W^\mu/\Gamma_{E_s}.$
If $l([w]) \neq l([w'])$ then $|\Delta\setminus I_{[w]}|$ and
$|\Delta\setminus I_{[w']}|$ differ at least by two.
Hence  $|I_{[w]}|$ and $|I_{[w']}|$ differ at least by two, so that by Theorem 3.3
$$Ext_{J(\Qp)}^1(v^J_{P_{I_{[w]}}},v^J_{P_{I_{[w']}}})=0.$$
Thus, the  expression (7) can be written as
\begin{eqnarray*} \nonumber  H^r(Y^{ad},\Upsilon)\otimes_{\Zl}\Ql\;  & \cong\; &  \bigoplus_{p \in \mathbb N}
  gr^p(H^r(Y^{ad},\Upsilon)\otimes_{\Zl}\Ql)  \\  \\  \nonumber 
 =  \bigoplus_{{[w]\in W^\mu/\Gamma_{E_s} \atop |\Delta\setminus I_{[w]}|=1 }\atop
  2l([w])=r}\!\!\!\!\!i^J_{P_{I_{[w]}}}\otimes\; ind_{[w]}(-l([w])) 
&\oplus& \!\!\!\!\!\! \bigoplus_{{[w]\in W^\mu/\Gamma_{E_s} \atop |\Delta\setminus
  I_{[w]}| >1} \atop 2l([w])=r}\!\!\!\!\!i^J_J\otimes\; ind_{[w]}(-l([w])) \\ \nonumber  &\oplus& \!\!\!\!\!\!\!\!\!
\bigoplus_{{[w] \in W^\mu/\Gamma_{E_s} \atop 2l([w]) + |\Delta\setminus I_{[w]}|-1=r} \atop p=|\Delta \setminus I_{[w]}|-1}\!\!\!\!\!\!\!\! v^J_{P_{I_{[w]}}}\otimes\; ind_{[w]}(-l([w])).
\end{eqnarray*}

\noindent {\bf Proof of Theorem 1.1:} The proof follows from the formula above together with the long exact sequence (5). \qed

\vspace{1cm}
\noindent Sascha Orlik, Universit\"at Leipzig, 
Mathematisches Institut, Augustusplatz 10/11, D-04109 Leipzig, Germany  

\noindent Email:orlik@mathematik.uni-leipzig.de

\end{document}